\documentclass[10pt]{amsart}

\usepackage{amsmath,amssymb,epsf, epic, eepic}
\usepackage{psfig}
\usepackage{amsfonts}
\usepackage{times}
\usepackage{graphicx}

\input xy
\xyoption{all}

\newtheorem{thm}{Theorem}[section]
\newtheorem{lem}[thm]{Lemma}
\newtheorem{cor}[thm]{Corollary}
\newtheorem{prop}[thm]{Proposition}

\theoremstyle{definition}

\newcommand{\bN}{\mathbb{N}}
\newcommand{\bZ}{\mathbb{Z}}

\newcommand{\bF}{\mathbb{F}}

\DeclareMathOperator{\Ker}{Ker}

\newcommand{\ra}{\rightarrow}

\newcommand{\ob}[1]{\overline{#1}}
\newcommand{\ul}[1]{\underline{#1}}

\newcommand{\spaces}{\;\;\;\;\;\;\;}
\newcommand{\ot}{\otimes}
\newcommand{\ol}{\overline}

\newcommand{\ei}{{\ul e}_i}

\newcommand{\FF}{\bF_2}

\newcommand{\kh}[3]{\mbox{KH}^{{#1},{#2}}({#3};\FF)}
\newcommand{\khn}[2]{\mbox{KH}_{n}^{{#1}}({#2};\FF)}
\newcommand{\khntri}[4]{\mbox{KH}_{n}^{{#1},{#2},{#3}}({#4};\FF)}
\newcommand{\chnu}[2]{\mbox{C}_{\ul n}^{{#1}}({#2})}
\newcommand{\sbn}[2]{\mbox{BN}^{{#1}} ({#2})}
\newcommand{\csbn}[2]{\mbox{BN}_{n}^{{#1}}({#2})}
\newcommand{\csbnu}[2]{\mbox{BN}_{\ul n}^{{#1}}({#2})}

\newcommand{\csbnbi}[3]{\mbox{BN}_{n}^{{#1},{#2}}({#3})}
\newcommand{\csbnbiu}[3]{\mbox{BN}_{\ul n}^{{#1},{#2}}({#3})}
\newcommand{\csbnbii}[4]{\mbox{BN}_{{#1}}^{{#2},{#3}}({#4})}

\newcommand{\chnm}[3]{C^{#1}_{#2}({#3})}
\newcommand{\hnm}[3]{H^{#1}_{#2}({#3})}
\newcommand{\chunum}[3]{C^{#1}_{\ul{#2}}({\ul{#3}})}
\newcommand{\chunumvar}[3]{C^{#1}_{{#2}}({{#3}})}
\newcommand{\hunum}[3]{H^{#1}_{\ul{#2}}({\ul{#3}})}
\newcommand{\hunumvar}[3]{H^{#1}_{{#2}}({{#3}})}
\newcommand{\khpolyn}{Kh_{n}}
\newcommand{\orient}{\theta}
\newcommand{\cano}{{\mathfrak s}_\theta}
\newcommand{\canoo}{{\mathfrak s}_{\ol{\theta}}}
\newcommand{\can}[1]{{\mathfrak s}_{{#1}}}

\newcommand{\framing}{f_{\!K}}
\newcommand{\symn}{{\mathfrak S}_n}
\newcommand{\symm}{{\mathfrak S}_m}
\newcommand{\symi}[1]{{\mathfrak S}_{#1}}
\newcommand{\symum}{{\mathfrak S}_{\ul m}}
\newcommand{\symun}{{\mathfrak S}_{\ul n}}
\newcommand{\sympart}[1]{({#1})^{\symn}}

\newcommand{\obd}{\ol{d}}

\newcommand{\uqsl}{U_q(\mathfrak{sl}_2)}
\newcommand{\lk}{\mbox{{\em lk}}}

\title{Bar-Natan's Khovanov homology for coloured links}
\author{Marco Mackaay}
\address{Departamento de Matem\'{a}tica\\ Universidade do Algarve\\ 
Campus de Gambelas\\ 8005-139 Faro\\ Portugal}
\email{mmackaay@ualg.pt}
\author{Paul Turner}
\address{School of Mathematical and Computer Sciences \\Heriot-Watt
  University\\ Edinburgh EH14 4AS\\Scotland}
\email{paul@ma.hw.ac.uk}
\begin{document}


\begin{abstract}
Using Bar-Natan's Khovanov homology we define a homology theory for
links whose components are labelled by irreducible representations of
$\uqsl$. We then compute this explicitly.
\end{abstract}

\maketitle


\vspace*{1cm}

\section{Introduction}
In \cite{khovanov2} Khovanov defined a link homology theory
categorifying the coloured Jones polynomial. He constructed a cochain
complex associated to an oriented framed link whose components are
labelled by irreducible representations of $\uqsl$ with the property
that the graded Euler characteristic of the homology of this complex 
is the coloured Jones polynomial. His key idea is to
interpret the formula
$$J_n(K)=\sum_{i=0}^{\lfloor\frac{n}{2}\rfloor}(-1)^i
\binom{n-i}{i}J(K^{n-2i}),$$
where $K^j$ is the $j$-cable of the knot $K$, as the 
Euler characteristic of a complex involving the link homology  of the
cablings $K^{n-2i}$, 
for $i=0,\ldots,\lfloor\frac{n}{2}\rfloor$. 

Cabling a knot or link immediately introduces an unmanagable number of
crossings from a computational point of view. Thus explicit
computations in coloured Khovanov theory are scarce. It is interesting
therefore to follow Khovanov's prescription for categorifying the
coloured Jones polynomial but in doing so replacing his original
theory with a simpler link homology theory. In this paper we do this
using the theory constructed by Bar-Natan in \cite{barnatan2}. To be
more precise we use the singly graded (filtered) theory which is
defined over $\FF$ obtained by setting the variable $H$ in
\cite{barnatan2} to be 1. In \cite{turner} this was computed
explicitly following the techniques of Lee \cite{lee}.

The resulting ``coloured'' invariants will not be particularly
interesting {\em per se} (cf computations of Lee's theory or Bar-Natan
theory where the homology only depends on the linking matrix of the
link), however as a testing ground for categorifying ``coloured''
invariants such a simple theory is invaluable. Also as Rasmussen's
remarkable paper \cite{rasmussen} shows, important topological
information can be extracted by considering the associated filtration
of a simple theory (in his case Lee's theory).

In Section \ref{sec:csbnt} we begin by recalling Khovanov's categorification of the coloured
Jones polynomial and follow this by the definition and calculation of Bar-Natan theory.  This categorification procedure can be carried out
for any link homology  and we proceed by using Bar-Natan
theory in place of Khovanov's theory
leading to the definition of what we call {\em coloured Bar-Natan theory}. We
calculate this explicitly for knots in Section \ref{sec:knots} where
the main results are Theorem \ref{thm:knots1} and Theorem \ref{thm:grading}. We extend this
to links in Theorem \ref{thm:links} of Section \ref{sec:links}. Finally we end with a
short section on an  alternative definition of the coloured theory
suggested by Khovanov.


\section{Defining coloured Bar-Natan theory}\label{sec:csbnt}

\subsection{Khovanov's categorification of the coloured Jones
  polynomial}\label{subsec:colkhov} Let $L$ be an oriented framed link
with each component labelled by an irreducible representation of
$\uqsl$. Such representations are parametrized by $\bN$, thus a
coloured link is a link with each component labelled or coloured by a
natural number. If this number is zero for a particular component we
can simply delete this component, so we can assume throughout that
$n\geq 1$. In this section we consider for the sake of simplicity an
oriented framed knot $K$ with colour $n$. Let $D$ be a diagram for $K$
whose blackboard framing corresponds to the given framing of $K$.

\begin{figure}[h]
\centerline{
\includegraphics{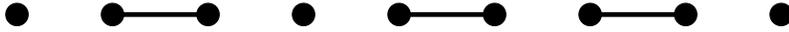}
}
\caption{A typical dot row}
\label{fig:dotrow} 
\end{figure}

A {\em dot-row} is a row of $n$ dots within which a number of
consecutive dots are paired. A typical dot-row can be seen in Figure
\ref{fig:dotrow}, where $n=9$.

Given a dot-row $s$ set $p(s) =$ number of pairs in $s$. Consider the
oriented graph whose vertices are dot rows and which has an arrow from  
$s$ to  $s^\prime$  if and only if all pairs in $s$ are pairs
  in $s^\prime$ and $p(s^\prime) = p(s) +1$.

\begin{figure}[h]
\centerline{
\includegraphics{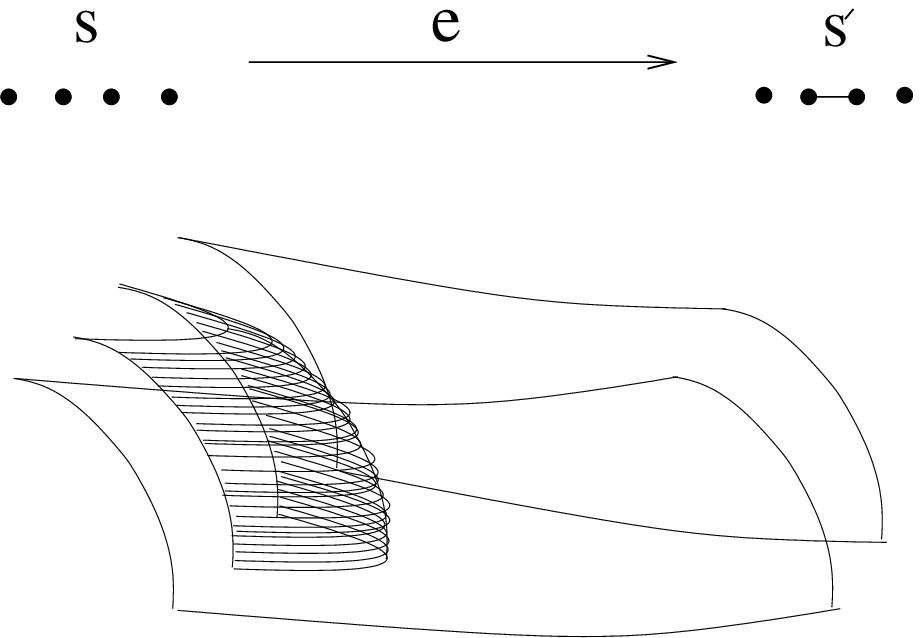}
}
\caption{}
\label{fig:edgecob} 
\end{figure}

Now attach to the vertex $s$ of this graph the link $D_s:= D^{n-2p(s)}$,
where here $D^k$ means the $k$-cable of the diagram $D$. Each single dot
in the dot-row corresponds to a cable-strand. Such a cable is oriented
such that adjacent cable-strands have opposite orientations. To an
edge $e$ with tail $s$ and head $s^\prime$ attach the cobordism $S_e$
which is the identity everywhere except at the two single dots in $s$
corresponding to the extra pair in $s^\prime$. For these two strands, the cobordism is
the annulus with two inputs and no outputs. This is illustrated in
Figure \ref{fig:edgecob}

The cochain complex for the coloured theory is now obtained by
applying Khovanov homology to this graph of links and
link-cobordisms. More precisely, set
\[
C^i = \bigoplus_s \kh ** {D_s}
\]
where the sum is over vertices $s$ such that $p(s) = i$. The
differential $d\colon C^i \ra C^{i+1}$ is defined as follows. Given a
dot row $s$ and $v\in \kh * * {D_s}$ set
\[
d(v) = \sum (S_e)_*(v) 
\]
where $(S_e)_*$ is the map in Khovanov homology induced by the
cobordism $S_e$.

Coloured Khovanov theory is then defined as the homology of this
complex, that is
\[
\khn i K = H_i(C^*,d).
\]

Each $C^k$ is bigraded via the bigrading in Khovanov homology
and this induced a trigrading on $\khn * K$. We
write
\[
\khn * K  = \bigoplus \khntri i j k K
\]
This defines a three variable polynomial
\[
\khpolyn (r ,t , q) = \sum_{i,j,k} r^i t^j q^k
\mbox{dim}_{\FF}(\khntri i j k K)
\]
from which one recovers the coloured Jones polynomial as
\[
J_n(K) = \khpolyn (-1,-1,q).
\]
In \cite{khovanov2} Khovanov shows that all the above is well
defined and can be extended to coloured links, yielding an invariant of
framed links.

Finally we remark that in the above $\kh **-$ may be replaced by any
Khovanov theory, by which we mean a theory satisfying the axioms given
by Bar-Natan in \cite{barnatan2} (see \cite{khovanov3} for a
discussion of the universal theory). For theories over rings other than
$\FF$ signs must be introduced and care must be taken with the sign
ambiguity introduced by the induced maps from cobordisms. As we will
be working over $\FF$ the above is enough for our purposes.

\subsection{ Bar-Natan Theory}\label{subsec:sbnt}
This is a singly graded link homology theory defined using the (ungraded)
Frobenius algebra $V=\FF \{1,x\}$ with multiplication $\tilde m$ given
by
\[
\tilde m (1,1) = 1 \spaces \tilde m (1, x) = x \spaces \tilde m (x, 1)=x
\spaces \tilde m (x, x) = x
\]
comultiplication $\tilde \Delta$
\[
\tilde\Delta(1) = 1\ot x + x\ot 1 + 1\ot 1 \spaces \tilde \Delta(x) = x\ot x
\] 
and unit and counit
\[
i(1) = 1 \spaces \epsilon (1) =0 \spaces \epsilon (x) =1.
\] 
Denote the resulting homology of a link $L$ by $\sbn * L$. We will refer
to this as {\em  Bar-Natan theory}. In fact Bar-Natan defines a bi-graded theory over $\FF[H]$ with deg$(H)=-2$. The theory we are considering is the assocaited filtered theory obtained by setting $H=1$. It is isomorphic to what is called stable Bar-Natan theory in \cite{turner}.

Bar-Natan theory was computed in \cite{turner} following the
techniques developed by Lee in \cite{lee}.

\begin{thm}
The dimension of $\sbn * L $ is $2^k$ where $k$ is the number of
components in $L$. 
\end{thm}

In fact explicit generators can be found as we now recall. There are
$2^k$ possible orientations of the diagram for $L$. Given an
orientation $\orient$ there is a canonical smoothing obtained by
smoothing all positive crossings to 0-smoothings and all negative
crossings to 1-smoothings. For this smoothing one can divide the
circles into two disjoint groups, Group 0 and Group 1 as follows. 
A circle belongs to Group 0 (Group 1) if it has the counter-clockwise 
orientation and is separated from infinity by an even (odd) number of circles 
or if it has the clockwise orientation and is separated from infinity by an 
odd (even) number of circles. Now consider the element in the chain complex
for $L$ defined by labelling each circle from Group 0 with $1+x$
and each circle from Group 1 with $x$, where we follow 
Bar-Natan's convention~\cite{barnatan1} for the ordering of the circles. 
It can be shown that this defines a cycle which we denote $\cano$ and refer 
to as the {\em canonical generator} from the orientation $\orient$. To prove the 
above theorem in fact one shows
\[
\sbn * L  \cong \FF \{ [\cano] \mid \orient \mbox{ is an orientation of $L$}\}
\]  

Furthermore, it is possible to determine the degree of the generators
in terms of linking numbers. Let $L_1, \ldots , L_k$ denote the
components of $L$. Recalling that $L$ is oriented from the start, if
we are given another orientation of $L$, say $\theta$, then we can
obtain $\theta$ by starting with the original orientation and then
reversing the orientation of a number of strands. Suppose that for the orientation
$\theta$ the subset $E\subset \{1,
2, \cdots , k \}$ indexes this set of  strands to be reversed. Let $\ol{E}=\{1,\ldots,n\}\backslash E$. The degree of the corresponding generator $[\cano]$ is then given by
\[
\mbox{deg}([\cano]) = 2 \times \sum_{l\in E,m\in \ol{E}} \lk(L_l,L_m)
\]
where  $\displaystyle{\lk(L_l,L_m)}$ is the linking number (for the original orientation) between component $L_l$
and $L_m$.

\subsection{Defining coloured  Bar-Natan theory}\label{subsec:sbn}
Coloured  Bar-Natan theory is the result of using $\sbn * -$
instead of $\kh ** -$ in the definitions in subsection
\ref{subsec:colkhov}. We now present the  general definition for links.

Let $L$ be an oriented framed link with $k$ components
and let $n_j\in \bN$ be the colour of the $j$'th component. It will be
convenient to assemble these into a vector $\ul n
= (n_1,n_2,\ldots , n_k)$.  Let $D$ be a
diagram for $L$ whose blackboard framing agrees with the given framing
of $L$ and let $D_{1}, D_{2}, \ldots , D_{k}$ be the components.

We will consider vectors of dot-rows $\ul s = (s_1, s_2, \ldots, s_k)$
where $s_j$ is a dot-row with $n_j$ dots. Let $\ul p (\ul s) =
(p(s_1), \ldots , p(s_k))$ and define $|\ul p (\ul s)| = p(s_1) +
\cdots + p(s_k)$ where as before $p(s_i)=$number of pairs in $s_i$.
Take the dot row vectors $\ul s$ as the vertices
of an oriented graph which has an arrow $\ul s \ra {\ul s}^\prime$ 
 if and only if all pairs in $\ul s$ are pairs in ${\ul s}^\prime$ and
 $|\ul p ({\ul s}^\prime)| = |\ul p ({\ul s})| + 1$.

 To a dot-row vector $\ul s$ we
attach the cable $D_{\ul s}= D^{\ul n - 2\ul p (\ul s)}$ where the
notation means that we take the $(n_j - 2p(s_j))$-cable of the
$j$'th component. Cables are oriented such that adjacent strands have
opposite orientation.  To an arrow $e\colon \ul s \ra {\ul s}^\prime$
we attach the cobordism $S_e$ as in subsection \ref{subsec:colkhov}
above i.e. we take the cobordism which is the identity everywhere
except at the two strands associated to the single dots in $\ul s$ corresponding to the extra
pair in ${\ul s}^\prime$. For these two strands, the cobordism is the
annulus with two inputs and no outputs.

To this we now apply
 Bar-Natan theory. Set
\[
\chnu i D = \bigoplus_{\ul s} \sbn * {D_{\ul s}}
\]
where the sum is over all $\ul s$ such that $|\ul p (\ul s)|= i$. 
This is a complex under the differential $d\colon \chnu i D \ra
\chnu {i+1} D$ defined for $v\in \sbn * {D_{\ul s}}$ by 
\[
d(v) = \sum (S_e)_*(v) 
\]
where $(S_e)_*$ is the map in  Bar-Natan theory induced by the
cobordism $S_e$ and the sum is over all edges $e$ with tail $\ul s$.

{\em Coloured  Bar-Natan theory} is then defined by 
\[
\csbnu i L = H_i(\chnu * D,d).
\]

Each $\chnu i D $ is graded which
induces a bigrading on coloured  Bar-Natan theory and we write
\[
\csbnu i L  = \bigoplus \csbnbiu i j L.
\]

It is easy to check that $\csbnbii {\ul 1} 0 * L \cong \sbn * L$, where $\ul 1 =
(1,1,\ldots , 1)$.


\section{Calculating coloured  Bar-Natan theory for knots}\label{sec:knots}
Let $K$ be a framed, oriented knot $K$ with framing $\framing\in \bZ$
and coloured by $n\in \bN$. Suppose this is presented via a
diagram $D$ with blackboard framing $\framing$.  Take a cross-section of the 
$n$-cable $D^n$ such that the original orientation of $D$ is upward and  
number the strands of $D^n$ with the numbers $1,\cdots, n$ from left to right. 
Recall that the $n$-cable is oriented such that adjacent strands have opposite 
orientation and we insist that the strand labelled 1 has the same 
orientation as $D$. 

Given any orientation $\orient$ of $D^n$, let $\orient (i)$ be
the orientation of the $i$'th cable. The symmetric group on $n$
letters $\symn$ acts on the set of orientations for $K^n$ as
follows. Given $\sigma\in\symn$ set
\[
\sigma\orient(i) = \orient (\sigma(i)).
\]
This induces an action of $\symn$ on the set of canonical generators
\[
\sigma \cano = \can {\sigma(\orient)}
\]
and hence an action on $\sbn * {K^n}$. This action is independent of the 
diagram we have chosen for $K$, because it commutes with the Reidemeister moves. 
This follows from Rasmussen's Prop. 2.3 in \cite{rasmussen}, which can easily be 
seen to hold for  Bar-Natan theory as well. In our case the proposition says 
that any Reidemeister move between two knot diagrams induces the linear isomorphism 
between the homologies of the two diagrams which maps any canonical generator to the 
canonical generator corresponding to the compatible orientation.   

Denote the vector space of symmetric elements under this action by
$\sympart {\sbn *{K^n}}$. That is
\[
\sympart {\sbn *{K^n}} = \{ \alpha \in \sbn *{K^n} \mid \sigma \alpha = \alpha \mbox{ for all }
\sigma \in \symn\}.
\]
It is clear that $\sympart
{\sbn *{K^n}}$ has a basis consisting of the elements $\sum [\cano]$ where the
sum is over a $\symn$-orbit of orientations $\orient$, from which it 
follows that $\sympart {\sbn *{K^n}}$ has dimension $n+1$.
Notice that using the diagram $D$ each $\symn$-orbit is determined soley by the number of
cable strands whose orientation agrees with the original (alternating)
orientation of $D^n$. 

In order to calculate the coloured  Bar-Natan theory we will
need the following lemma. Let $e$ be an edge of the graph of dot-rows
connecting $s$ to $s^\prime$ and let $S_e$ be the associated cobordism
as defined in section \ref{subsec:sbn}. Suppose that strands $l$ and
$l+1$ are contracted. An orientation $\orient$ of $L_s$ induces an
orientation $\orient^\prime$ on $L_{s^\prime}$, simply by removing
strands $l$ and $l+1$.

\begin{lem}\label{lem:cangen}
\[
(S_e)_* ([\cano]) = \begin{cases}
[\can {\orient^\prime}] & \mbox{ if $l$ has opposite orientation to
  $l+1$}\\
0 &  \mbox{ if $l$ and $l+1$ have the same orientation}
		    \end{cases}
\]
\end{lem}

\begin{proof}
This lemma is essentially a variation mod 2 on a result by 
Rasmussen. The cobordism $S_e$ can be presented by a movie starting
with a 1-handle move which fuses two strands into one, followed by a
sequence of type two Reidemeister moves, and finally ending with the
removal of a circle. 

The first of these is given at the chain level by multiplication or
comultiplication. If the orientations of the two strands agree, then
they belong to different circles and the cobordism induces a
multiplication. Based on Rasmussen's analysis (Lem. 2.4
in~\cite{rasmussen}) we see that the canonical generator will label
one of these strands with $1+x$ and the other with $x$, thus their
fusion will produce zero as required. If the orientations are opposite
there are two possible cases. In the first case the two strands belong
to different circles and the cobordism induces a multiplication. The
label will be the same on the two and fusion will give that same label
to the single fused circle. In the second case the two strands belong
to the same circle, with one lable of course, and the cobordism
induces a comultiplication which produces two circles with that same
lable. In both cases we see that a canonical generator is mapped to a
canonical generator.

As already remarked above, by arguments similar to Rasmussen's (proof
of Prop. 2.3 in \cite{rasmussen}) it can be seen that Reidemeister
two moves take a canonical generator to a canonical generator.

Finally the removal of a circle is given by the counit which takes the
value 1 on both $1+x$ and $x$. 
\end{proof}

\begin{cor}\label{cor:perm}
Let $\sigma$ be the permuation
switching strands $l$ and $l+1$.  Then 
\[
 (S_e)_* ([\cano])= (S_e)_*
([\can {\sigma \orient}])
\]
\end{cor}

We now
calculate the coloured  Bar-Natan theory in degree zero for a
coloured knot. 

\begin{prop} \label{prop:sym}
Let $K$ be a framed, oriented knot with colouring $n\in\bN$. Then
\[
\csbn 0 K = \sympart {\sbn *{K^n}} 
\]
\end{prop}

\begin{proof}
We must show $\Ker (d_0) = \sympart {\sbn *{D^n}}$ where $d_0$ is the
differential $d_0\colon C^0_n(D) \ra C^1_n(D)$.  Notice that $d_0$ is the sum of
$n-1$ contractions. Let $d_{0,l}$ be the induced map in  Bar-Natan
theory corresponding to the cobordism contracting strands $l$ and
$l+1$. An element $v$ is in the kernel of $d_0$ if and only if it is
in the kernel of $d_{0,l}$ for $l=1, \cdots, n-1$.

Now let $v = \sum [\cano] \in \sympart {\sbn *{K^n}}$. We claim that
$d_{0,l}(v) = 0$. Firstly note that we can split $v$ into a sum over 
orientations where strands $l$ and $l+1$ are oriented the same way and
orientations where they are oriented differently. 
\[
v = \sum_{\scriptstyle\text{same}} [\cano] + \sum_{\scriptstyle\text{diff}} 
[\cano]
\] 
Letting $\sigma$ be the transposition switching $l$ and $l+1$ the sum on
the right can be replaced by a sum of elements $[\cano] + [\sigma
  \cano]$ where the sum is over orientations for which strand $l$ and
$l+1$ are oriented as for the
original orientation of $K^n$.
\[
v = \sum_{\scriptstyle\text{same}} [\cano] + \sum_{\scriptstyle\text{orig}}( [\cano] + [\sigma \cano])
\]
Thus 
\begin{eqnarray*}
d_{0,l} (v) & = & d_{0,l} (\sum_{\scriptstyle\text{same}} [\cano] + \sum_{\scriptstyle\text{orig}}( [\cano] +
[\sigma \cano]))\\
 & = & \sum_{\scriptstyle\text{same}} d_{0,l} ([\cano]) + \sum_{\scriptstyle\text{orig}} (d_{0,l} ([\cano]) +
d_{0,l} ([\sigma \cano])) \\
 & = & \sum_{\scriptstyle\text{orig}} d_{0,l} ([\cano]) +
d_{0,l} ([\sigma \cano]) \;\;\;\;\mbox{by Lemma \ref{lem:cangen}}\\
 & = & 0 \;\;\;\;\mbox{mod 2} \;\;\;\;\;\;\;\;\mbox{by Corollary \ref{cor:perm}}
\end{eqnarray*}
Since this is true for $l=1,\cdots, n-1$ we have shown $d_0(v) =0$ and
hence  $\sympart {\sbn *{K^n}} \subseteq \Ker (d_0)$.

Now suppose that $v = \sum \lambda_\orient [\cano]$ is in the kernel
of $d_0$ and let $l\in \{1, \cdots , n-1\}$. As above let $\sigma$ be the transposition 
switching strands $l$ and $l+1$. Suppose for an orientation
$\orient$ we have $\lambda_\orient \neq 0$. If strands $l$ and $l+1$
have the same orientation in $\theta$ then $\sigma [\cano] = [\cano]$. If strands
$l$ and $l+1$ have different orientations, then using the fact that
$v\in\Ker (d_{0,l})$, Lemma \ref{lem:cangen} implies that we see we must also have
$\lambda_{\sigma \orient} \neq 0$. This gives us that $\sigma v = v$. 
This is true for all transpositions 
$\sigma= (l,l+1)$ for $l=1,\cdots , n-1$ and since $\symn$
in generated by such transpositions this shows that $\tau v = v$ for all
$\tau \in \symn$. We conclude that $v\in \sympart {B^*(K^n)}$ and
hence $\Ker (d_0) \subseteq \sympart {\sbn *{K^n}}$, which
finishes the proof.
\end{proof}

To calculate the coloured  Bar-Natan theory in higher degrees we
need to first introduce a new family of complexes. For $m=0,1,\ldots ,
n$ define $\chnm * n m $ as follows, where we surpress the diagram  $D$
from the notation. Consider the graph of dot-rows with
$m$ dots. Now attach to the vertex $s$ of this graph the cable $D_s:=
D^{n-2p(s)}$. The first $(m-2p(s))$-strands of this cable are
associated with the $m-2p(s)$ dots. The remaining $n-m$ cable strands
are not associated to any dots. To edges we associate annulus  cobordisms as
before insisting that cobordisms  are the identity on the last $n-m$ ``free'' cables.
Now set
\[
\chnm i n m = \bigoplus_{s} \sbn * {D_s}
\]
where the sum is over dot-rows $s$ such that $p(s) = i$. As above the
arrows in the dot row graph give rise to a differential. Let $\hnm * n m$ denote
the homology of this complex.

Notice that $\chnm * {n-2}{m-2}[-1]$ can be identified with the
subcomplex of $\chnm * n m $ consisting of dot-rows with the last two
dots always forming a pair. Here $[-1]$ means a downward shift of one
degree. The quotient complex is seen to be $\chnm * n
{m-1}$ and so there are short exact sequences
\[
\xymatrix{ 0 \ar[r] & \chnm * {n-2} {m-2} [-1] \ar[r] & \chnm * n m  \ar[r] &
  \chnm  * n {m-1} \ar[r] & 0}.
\]

\begin{figure}[h]
\centerline{
\includegraphics{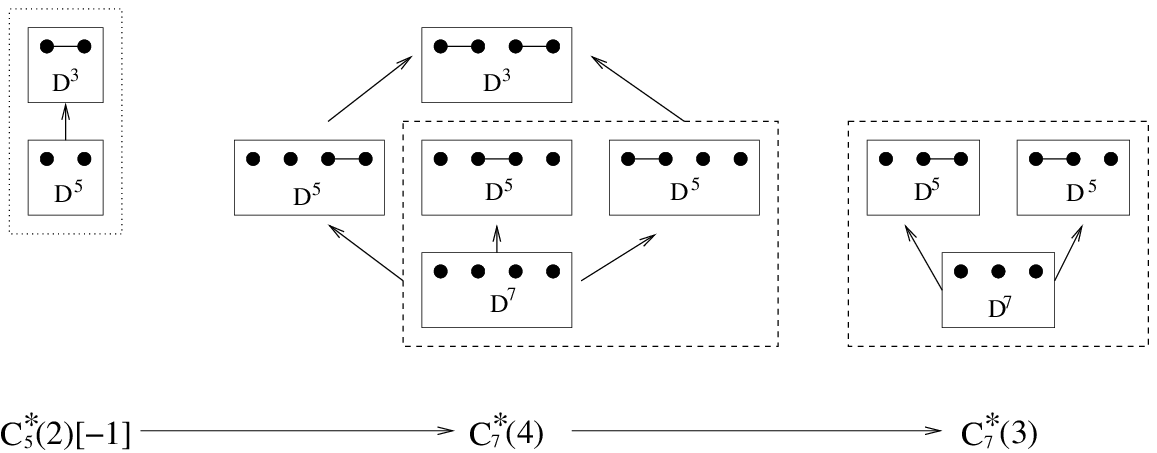}
}
\caption{}
\label{fig:ses} 
\end{figure}

This is illustrated for $n=7$ and $m=4$ in Figure \ref{fig:ses}.

Such a short exact sequence gives rise to a long exact sequence as follows.
{\small
\[
\xymatrix{\cdots \ar[r] & \hnm {i-1} {n-2} {m-2} \ar[r] & \hnm i n m
  \ar[r] & \hnm i n {m-1} \ar[r]^-{\delta} & \hnm i {n-2}{m-2} \ar[r] & \cdots }
\]
}
In the next lemma we gather some useful results about
the homology of these complexes. In part (3) we are using the fact
that $\symm$ includes in $\symn$ as permutations of the first $m$ members of
$\{1, \ldots , n\}$ and so $\symm$ acts on $\sbn * {K^n}$. 

\begin{lem}\label{lem:hnm}
\begin{eqnarray}
\hnm * n n & \cong & \csbn * K \\
 \hnm i n 0 &\cong &\hnm i n 1 \cong \begin{cases}
\sbn * {K^n} & i=0\\ 0 & \mbox{otherwise} 	 \end{cases}\\
\hnm 0 n m & \cong & (\sbn * {K^n})^{\symm} \\ 
\mbox{dim} (\hnm 0 n m) & = & 2^{n-m}(m+1)
\end{eqnarray}
\end{lem}

\begin{proof}
Examining the definition above shows that $\chnm * n n$ is the complex
defining coloured  Bar-Natan theory and also that 
\[
\chnm i n 0 \cong \chnm i n 1 \cong \begin{cases}
\sbn * {K^n} & i=0\\ 0 & \mbox{otherwise} 	 \end{cases}\]
from which (1) and (2) follow immediately. Part (3) is essentially
a corollary to Proposition \ref{prop:sym} and Part (4) follows from
part (3).
\end{proof}

\begin{prop}\label{prop:knotindbase}
$\hnm 1 n m = 0$ for all $n\in\bN$ and $m=0,1,\ldots , n$.
\end{prop}

\begin{proof}
The proof is by induction on $m$. By Lemma \ref{lem:hnm} (2) we have $\hnm 1
n 0 = \hnm 1 n 1 = 0$ for all $n$ which starts the induction. Now
suppose that $\hnm 1 n m =0$ for $m=0,1,\ldots , M$ and $n\geq
M$. We claim $\hnm 1 n {M+1} = 0 $ for $n\geq M+1$. There is a piece
of long exact sequence as follows.
{\small
\[
\xymatrix{0 \ar[r] & \hnm 0 n {M+1} \ar[r] & \hnm 0 n {M} \ar[r] & \hnm 0 {n-2}{M-1} \ar[r] &
\hnm 1 n {M+1} \ar[r] & \hnm 1 n {M} }
\]
}
The group on the right is trivial by inductive hypothesis. Using Lemma
\ref{lem:hnm} (4) we see that the first three groups in the sequence
have dimensions $2^{n-M-1}(M+2)$, $2^{n-M}(M+1)$ and $2^{n-M-1}M$
respectively. Since the alternating sum of dimensions in an exact
sequence is zero this gives
\begin{eqnarray*}
\mbox{dim}(\hnm 1 n {M+1})  & = &   2^{n-M-1}(M+2) - 2^{n-M}(M+1) +
2^{n-M-1}M \\ & = &  2^{n-M-1}(2M+2) - 2^{n-M}(M+1) = 0.
\end{eqnarray*}
Hence $\hnm 1 n {M+1} = 0$ as required.
\end{proof}

\begin{prop}\label{prop:knothigher}
$\hnm i n m = 0$ for all $i\geq 1$ and for all $n\in\bN$ and $m=0,1,\ldots , n$.
\end{prop}

\begin{proof}
The proof is by induction on $i$. Proposition \ref{prop:knotindbase} starts the
induction with $i=1$. Now suppose the result holds for $1\leq i <
j$. We claim that the result holds for $i=j$, i.e. $\hnm j n m = 0$ for
all $n$ and $m=0,1,\ldots , n$. To show this we do  a second induction
this time on $m$ where the base
case is again provided by Lemma \ref{lem:hnm} (2). Now suppose $\hnm
j n m = 0$ for $m=0,1,\ldots , M$. Then we have a piece of long exact
sequence
\[
\xymatrix{ \hnm {j-1} {n-2}{M-1} \ar[r] & \hnm j n {M+1} \ar[r] & \hnm j n {M} }.
\]
The group on the right is trivial by the induction on $m$ and the
group on the left is trivial by the induction on $i$, hence the middle
group is trivial too.
\end{proof}

In particular when $m=n$ we have $\csbn i K = \hnm i n n = 0$ for
$i\geq 1$. Combining this with Proposition \ref{prop:sym} we have the
following theorem.

\begin{thm}\label{thm:knots1}
\[
\csbn i K = \begin{cases}
\sympart {\sbn *{K^n}} & i= 0\\
0 &  \mbox{otherwise}
		      \end{cases}
\]
\end{thm}

Theorem~\ref{thm:knots1} is similar to Khovanov's Theorem 1 in 
\cite{khovanov2}. However, we note that our proof, although similar in spirit, is 
more complicated 
due to the fact that $\sbn *{K^n}\not\cong \sbn*{K}^{\otimes n}$.

\begin{cor}
\[
\mbox{dim}(\csbn * K) = n+1.
\]
\end{cor}

Recall that in each degree there is a second grading. For a complete
calculation it remains for
us to describe this for the zero'th homology group. Note that if
$\cano\in \sbn i{K^n}$ then for $\sigma\in \symn$ we have $\can
{\sigma\orient} \in \sbn i{K^n}$, because the number of negative crossings 
in $K^n$ is the same for both orientations, which shows us that the grading on
$\sbn * {K^n}$ induces a well defined internal grading on $ \csbn 0 K$. We
write $\csbn 0 K = \bigoplus \csbnbi 0 i K$.

\begin{thm}\label{thm:grading}
If $n$ is even then
\[
\csbnbi 0 i K = \begin{cases}
\FF & i = 0 \\
\FF \oplus \FF & i = -2k^2 f_K \;\;\;\;\; k = 1, \ldots, \frac{n}2\\
0 & \mbox{otherwise}
		      \end{cases}
\]

If $n$ is odd then
\[
\csbnbi 0 i K = \begin{cases}
\FF \oplus \FF & i = -2k(k+1) f_K \;\;\;\;\; k = 0, \ldots, \frac{n-1}2\\
0 & \mbox{otherwise}
		      \end{cases}
\]
\end{thm}
\begin{proof}
Suppose that $n$ is even. In this case the orientation given to the
$n$-cable of $K$ has $\frac{n}2$ strands agreeing with the given
orientation of $K$ and $\frac n 2$ strands with the opposite
orientation. Recall that an $\symn$-orbit is determined by the number
of strands agreeing with the original orientation. Moreover, corollary~\ref{cor:perm} 
implies that any element in a given orbit has the same internal degree. Note 
that $S_e$ preserves the internal degrees. If the orbit of $\orient$ has more than 
one element, then there exists a transposition $\sigma$ such that $\orient$ and 
$\sigma \orient$ are different. By corollary~\ref{cor:perm} we know that 
$\cano$ and $\can{\sigma \orient}$ are mapped to the same non-zero element, so they 
have to have the same degree. 

The given orientation of $K^n$ has $\frac n 2$ strands agreeing with
$K$ and using the description above in Subsection \ref{subsec:sbnt} it
is immediate that the canonical generator in $\sbn *{K^n}$
corresponding to this orientation has degree zero (the set $E$ in the
theorem is the empty set). Thus the class in $\csbn 0 K$ determined by
the orbit of this has internal grading zero.

 Let $D_1, \ldots , D_n$
be the strands of the $n$-cable diagram and note that $\lk (D_{2i}, D_l)
= \lk (D_2,D_l)$ and moreover that this is equal to $\framing$ if $l$
is even and $-\framing$ if $l$ is odd.

Now let $E=\{2,4, \ldots , 2k\}$ for $k = 1,2, \ldots , \frac n 2$ and
consider the orientation of $K^n$ obtained by reversing the
orientation of strands numbered by elements of $E$. Such an
orientation has $\frac n 2 + k$ strands agreeing with the orientation
of $K$. We see that the degree of the canonical generator defined by
this orientation is given by
\begin{eqnarray*}
2 \sum_{j\in E, l\in \ol E} \lk (D_j, D_l) & = & 2 (\sum_{l\in \ol E} \lk
(D_2,D_l) + \cdots + \lk (D_{2k},D_l)) \\
& = & 2k\sum_{l\in \ol E} \lk (D_2,D_l)\\
 & = &  2k(\sum_{l\in \{2k+2, \ldots , n\}} \lk (D_2,D_l)+ \sum_{l\,
\scriptstyle\text{odd}} \lk (D_2,D_l))\\
 & = &  2k(\sum_{l\in \{2k+2, \ldots , n\}} \framing + \sum_{l\,
  \scriptstyle\text{odd}} (-\framing) )\\
 & = & -2k^2\framing
\end{eqnarray*}

Thus the $\symn$-orbit containing this orientation has grading $-2k^2\framing
$. This accounts for the $\symn$-orbits with $> \frac n 2$ strands
agreeing with $K$.

To account for the orbits with $<\frac n 2 $ strands agreeing with
$K$, let $E=\{1,3, \ldots , 2k-1\}$ for $k = 1,2, \ldots , \frac n
2$ and consider the orientation of $K^n$ obtained by reversing the 
orientation of strands in $E$. This orientation has
$\frac n 2 - k$ strands agreeing with $K$. By a similar computation to
the above we see that  the degree of the canonical
generator defined by this orientation is given by $-2k^2\framing$
thus, combining with the above, giving a total of two generators in
this degree. 

This finishes the proof for $n$ even.

For $n$ odd the argument is similar with the difference that the given
orientation of $K^n$ and the orientation from $E=\{1\}$ both give classes in
degree zero. The orientations from $E=\{2,4, \ldots , 2k\}$ and $E =
\{ 1, 3, \ldots , 2k+1\}$ for $k=1,2, \ldots \frac{n-1}2$ both give
rise to classes in degree $-2k(k+1)\framing$. A slight re-arrangment
of indexing gives the statement of the theorem for $n$ odd.
\end{proof}

Recall that each generator $\alpha$ is of the form $\sum_{\theta\in X}
[\cano]$ where $X$ is a $\symn$-orbit of orientations. It is easy to
see that $\ob \alpha := \sum_{\theta\in X} [\canoo]$ is also a
generator, because $\sigma \ol \orient=\ol{\sigma\orient}$, for any $\sigma\in\symn$. 
Moreover, by inspecting
the above proof one can see that both $\alpha$ and $\ob \alpha$ have
the same grading. Thus for each $\FF \oplus \FF$ in Theorem
\ref{thm:grading} if $\alpha$ generates one copy of $\FF$ then
$\ob \alpha$ generates the other. Note that if $n$ is even and $\alpha$ is 
the generator of degree 0, then $\ol \alpha=\alpha$.


\section{Calculating coloured  Bar-Natan theory for
  links}\label{sec:links}

Let $L$ be an oriented, framed link with $k$ components and let $n_j\in\bN$
be the colouring of the $j$'th component. Let $D$ be a diagram for $L$
whose blackboard framing is the given framing of $L$. 

As above for the case of knots we need a new family of complexes in order to
calculate the coloured  Bar-Natan theory of $L$. For each
$j=1,\ldots, k$ choose $m_j\in \bN$ and set $\ul m =
(m_1,\ldots , m_k)$. We say that $\ul m$ is {\em allowable} if $1\leq m_j \leq n_j$ for
$j=1,\ldots , k$. We define the {\em length} of $\ul m$ by $|\ul m
| = m_1+ \cdots +m_k$.

Now consider the oriented graph whose vertices are vectors of dot rows
$\ul s = (s_1, \ldots , s_k)$ where $s_j$ is a dot row of $m_j$
dots. Arrows in the graph are defined using the same construction as in subsection
\ref{subsec:sbn} coordinatewise. Now attach $D_{\ul s} = D^{\ul n - 2 \ul p (\ul s)}$
to the vertex $\ul s $, where recall that $\ul p (\ul s) = (p(s_1),
\ldots , p(s_k))$.  For the $j$'th component, the first $m_j-2p(s_j)$
strands of $D^{n_j - 2 p ( s_j)}$ are identified with the
$m_j-2p(s_j)$ dots, and the remaining $n_j-m_j$ strands are ``free''.
As usual attach cobordisms to edges.  Now set
\[
\chunum i n m = \bigoplus_{\ul s} \sbn * {D_{\ul s}}
\]
where the sum is over vectors of dot-rows $\ul s$ such that $|\ul p(\ul
s)| = p(s_1) +
\cdots + p(s_k) = i$. As above the
arrows in the graph give rise to a differential. Let $\hunum * n m$ denote
the homology of this complex.

Let ${\ul e}_i = (0,\ldots , 0 ,1,0,\ldots ,0)$ with a $1$ in the
$i$'th place. Since the arrows in the graph are defined coordinatewise, there are 
short exact sequences 
\[
\xymatrix{ 0 \ar[r] & \chunumvar * {\ul n-2\ei} {\ul m-2\ei} [-1] \ar[r] & \chunum * n m  \ar[r] &
  \chunumvar  * {\ul n} {\ul m-\ei} \ar[r] & 0}.
\]
giving rise to long exact sequences
{\small
\[
\xymatrix{\cdots \ar[r] & \hunumvar {i-1} {\ul n-2\ei} {\ul m-2\ei} \ar[r] & \hunum i n m
  \ar[r] & \hunumvar i  {\ul n} {\ul m-\ei} \ar[r]^{\delta} &
  \hunumvar i  {\ul n-2\ei} {\ul m-2\ei} \ar[r] & \cdots } 
\]
}

The symmetric group $\symi {n_j}$ acts on the orientations of the
$n_j$-cable of the $j$'th component and so $\symi {\ul n} = \symi
{n_1}\times \cdots \times \symi {n_k}$ acts on the set of orientations
of $L^{\ul n}$. As before this induces an action on $\sbn * {L^{\ul
    n}}$. Letting $\symi {\ul m} = \symi {m_1} \times \cdots \times
\symi {m_k} \subset  \symi {n_1}\times \cdots \times \symi {n_k} =
\symi {\ul n}$ we see that $\symi {\ul m}$ acts on $\sbn * {L^{\ul
    n}}$ too. We will write $ (\sbn * {L^{\ul n}})^{\symum} $ to
denote the space of symmetric elements under this action.

 There is an analogue of Lemma \ref{lem:hnm} for
links.

\begin{lem}\label{lem:hunum}$~$\\[-4mm]
\begin{enumerate}
\item $\hunum * n n  \cong  \csbnu * L $
\item If  $\ul m$ is such that $0\leq m_j \leq 1$ for $j=1,\ldots , k$
  then 
\[
 \hunum i n m \cong \begin{cases}
\sbn * {L^{\ul n}} & i=0\\ 0 & \mbox{otherwise} 	 \end{cases}
\]
\item If $\ul m$ is allowable then 
\[
\hunum 0 n m \cong  (\sbn * {L^{\ul n}})^{\symum} 
\] 
\item \[ \mbox{dim}(\hunum 0 n m)  =  \prod_{j=1}^k2^{n_j-m_j}(m_j +1)\]
\end{enumerate}
\end{lem}

\begin{proof}
The only non-trivial thing to check is part (3) which is an easy
generalisation  of Prop \ref{prop:sym}.
\end{proof}

The following results are the analogues of Propositions~\ref{prop:knotindbase} and  
\ref{prop:knothigher} and 
Theorem~\ref{thm:knots1}.

\begin{prop}\label{prop:linksindbase}
$\hunum 1 n m = 0$ for all $\ul n$ and all allowable $\ul m$.
\end{prop}

\begin{proof}
The proof is by induction on the length of $\ul m$. The induction is
started by Lemma \ref{lem:hunum} (2). Suppose that $\hunum 1 n m = 0
$ for all allowable $\ul m$ with $|\ul m|<M$. Now let $\ul m$ be such
that $|\ul m|=M$. Using Lemma \ref{lem:hunum} (2) we can suppose that
for some $i$ we have $m_i\geq 2$. There is a piece of long exact
sequence
{\small
\[
\xymatrix{0 \ar[r] & \hunumvar 0 {\ul n} {\ul m} \ar[r] & \hunumvar 0
  {\ul n} {\ul m - \ei} \ar[r] & \hunumvar 0 {\ul n-2\ei}{\ul m -2 \ei} \ar[r] &
\hunum 1 n m \ar[r] & \hunumvar 1 {\ul n} {\ul m -\ei} }
\]
}
Since $|\ul m - \ei | = M-1$ the group on the right is zero by the
inductive hypothesis. Using Lemma \ref{lem:hunum} (4) the first three
groups have dimensions as follows.
\begin{eqnarray*}
\mbox{dim}(\hunum 0 n m)  & = & \prod_{j}2^{n_j-m_j}(m_j +1)\\
\mbox{dim}(\hunumvar 0 {\ul n} {\ul m - \ei})  & = &  2^{n_i-m_i +
  1}m_i\prod_{j\neq i}2^{n_j-m_j}(m_j +1)\\ 
\mbox{dim}(\hunumvar 0 {\ul n-2\ei} {\ul m-2\ei})  & = &
2^{n_i-m_i}(m_i -1)\prod_{j\neq i}2^{n_j-m_j}(m_j +1)
\end{eqnarray*}
Using the fact that the alternating sum of dimensions in an exact
sequence is zero we get
\begin{eqnarray*}
\mbox{dim}(\hunum 1 n m )
 & = &  \Bigl[ 2^{n_i - m_i}(m_i +1) - 2^{n_i-m_i +
  1}m_i + 2^{n_i-m_i}(m_i -1)\Bigr] \prod_{j\neq i}2^{n_j-m_j}(m_j +1)\\
 & = & 2^{n_i - m_i}[m_i +1 - 2m_i
  + m_i -1]\prod_{j\neq i}2^{n_j-m_j}(m_j +1) = 0.
\end{eqnarray*}
This shows $\hunum 1 n m = 0 $ as required.
\end{proof}

\begin{prop}\label{prop:linkshigher}
$\hunum i n m = 0$ for all $i\geq 1$ and for all $\ul n$ and allowable $\ul m$.
\end{prop}

\begin{proof}
We use induction on $i$ which is started by Proposition
\ref{prop:linksindbase}. Suppose $\hunum i n m = 0$ for $1\leq i < j$. We 
now claim that $\hunum j n m = 0$. To show this we use a second
induction this time on the length of $\ul m$ which is started by Lemma
\ref{lem:hunum} (2). Suppose $\hunum j n m = 0$ for $|\ul m|<M$. We 
claim that for $\ul m$ satisfying $|\ul m | = M$ we have $\hunum j n m =
0$. There is a piece of long exact sequence
\[
\xymatrix{\hunumvar {j-1} {\ul n - 2\ei} {\ul m -2\ei}
  \ar[r] & \hunumvar j  {\ul n} {\ul m} \ar[r] &
  \hunumvar j  {\ul n} {\ul m-\ei} }. 
\] 
The group on the right is zero by the second induction since $|\ul m
-\ei|<M$ and the group on the left is zero by the first induction,
thus the middle group is zero too.
\end{proof}

\begin{thm}\label{thm:links}
\[
\csbnu i L = \begin{cases}
(\sbn *{L^{\ul n}})^{\symun} & i= 0\\
0 &  \mbox{otherwise}
		      \end{cases}
\]
\end{thm}
\begin{proof} From Lemma \ref{lem:hunum} (1) and (3) we observe that 
$\csbnu 0 L \cong \hunum 0  n
n \cong (\sbn * {L^{\ul n}})^{\symun}$. 
Proposition~\ref{prop:linkshigher} shows that $\csbnu i L \cong \hunum i n
n = 0$, for all $i\geq 1$.
\end{proof}

\begin{cor}
\[
\mbox{dim}(\csbnu * L) = \prod_{j=1}^k (n_j+1).
\]
\end{cor}

We now consider the generators of $\csbnu * L $. For each $i$ let $E_i$ be a subset of $\{1,\ldots ,
n_i\}$ and define $\ob E_i = \{1,\ldots ,n_i\} \setminus E_i$. We
denote the collection $\{ E_i\}_{i=1,\ldots , k} $ by $\ul E$. Given
$\ul E$ define 
\begin{eqnarray*}
\lambda_{i,l} & = &
\#\{(j,m)\in  E_i\times   \ob E_l \;|\; j \mbox{ and } m
\,\mbox{have the same parity}\} \\  
& & \quad \quad \quad  \quad- \#\{(j,m)\in E_i\times \ob E_l \;|\; j \mbox{ and } m
\,\mbox{have different parities}\}\\
& = &  \sum_{j\in E_i}\sum_{m\in \ob E_l} (-1)^{j+m}
\end{eqnarray*}

The collection $\ul E$ is said to be {\em admissible} if each $E_i$ is
either empty or of one of the following forms:
\begin{eqnarray*}
E_i & = & \{ 2,4, \ldots , 2p\} \mbox{ for some } p \in \{ 1, \ldots ,
  \lfloor\frac{n_i}2 \rfloor \}\\
E_i & = & \{ 1,3, \ldots , 2p-1\} \mbox{ for some } p \in \{1, \ldots ,
  \lceil \frac{n_i}2 \rceil \}
\end{eqnarray*}

Each admissible $\ul E$ determines an orientation of $L^{\ul n}$ by
reversing the orientation of the strands indexed by $E_i$ in the $n_i$-cable of the
$i$'th component. As explained in subsection \ref{subsec:sbnt}  this
defines a cocycle $s_{\ul E}$ and hence a class $[s_{\ul E}]\in \sbn *
{L^{\ul n}}$. The sum of elements in the $\symun$-orbit of $
[s_{\ul E}]$ is a generator of $\csbnu 0 L \cong (\sbn * {L^{\ul
    n}})^{\symun}$. All generators arise in this way and there is a
  one-to-one correspondence between the admissbile $\ul E$ and
  generators of $\csbnu 0 L$.

Recall that $\csbnu * L$ is in fact bi-graded, the second grading
being inherited from the grading on  Bar-Natan theory. The following 
theorem is the analogue of Theorem~\ref{thm:grading} for knots.

\begin{thm}
Let $\ul E$ be an admissible collection. The generator of $\csbnbiu 0
* L$ corresponding to $\ul E$ has grading
\[
2 \sum_{i=1}^k\sum_{l=1}^k \lambda_{i,l} \lk (L_i,L_l).
\]
\end{thm}

\begin{proof}
It suffices to work out the degree of $[s_{\ul E}]$, because the generators in one 
orbit all have the same internal degree as explained in the proof of 
Theorem~\ref{thm:grading}. Let
\[
E = \{ (i,j) \; | \; i=1,\ldots , k \mbox{ and } j\in E_i\}
\]
and set
\[
\ob E = \{ (i,j) \; | \; i=1,\ldots , k \mbox{ and } j=1,\ldots ,
n_i\} \setminus E.
\]
By the discussion of  Bar-Natan theory the degree of $[s_{\ul
  E}]$ is
\[
2\sum_{\substack{(i,j)\in E \\ (l,m)\in \ob E}} \lk (D_{ij},D_{lm}).
\]
where $D_{ij}$ is the $j$'th strand of the $n_i$-cable of the $i$'th
component of the diagram $D$. Recalling that the strands of each cable
are oriented such that strand 1 has the given orientation and that
adjacent strands have opposite orientations we see that $\lk (D_{ij},
D_{lm}) = (-1)^{j+m}\lk (L_i, L_l)$.  Thus we have
\begin{align*}
\mbox{degree}( &[s_{\ul E}])  =  2\sum_{\substack{(i,j)\in
    E\\(l,m)\in \ob E}} \lk  (D_{ij},D_{lm})
  =  2\sum_{i=1}^k\sum_{l=1}^k\sum_{\substack{j\in E_i \\ m\in \ob
    E_l}} \lk (D_{ij},D_{lm})\\
 & =  2\sum_{i=1}^k\sum_{l=1}^k\sum_{\substack{j\in E_i\\m\in \ob
    E_l}} (-1)^{j+m}  \lk (L_i,L_l)
 =   2\sum_{i=1}^k\sum_{l=1}^k\lambda_{i,l}\lk (L_i,L_l)
\end{align*}
\end{proof}


\section{Khovanov's alternative definition}\label{sec:kc}
There is an alternative way of defining a complex which categorifies
the coloured Jones polynomial. This method, suggested by Khovanov,
involves reversing the direction of the arrows in the graph of dot-rows.
Khovanov conjectures that working over a field of characteristic zero
with the Khovanov homology the two complexes give isomorphic homology
groups. In this section we prove the analogous statement for the coloured 
 Bar-Natan theory. We note that our proof is valid in Khovanov's 
original setting too.

Let $(C^*,\obd)$ be the chain complex obtained from the reversed
 graph. The reversed graph has cables attached to vertices as before
 and cobordisms attached to edges again except this time with the
 opposite orientation. This way we get a chain complex (rather than a
 cochain complex) with differential $\obd_{i} \colon C^{i} \ra
 C^{i-1}$.

We can identify $\obd$ in terms of $d$.  Given any link $L$ there is a
non-degenerate inner product on $\sbn *L$. This inner product is
induced by the cylinder cobordism from the link $L \sqcup \ol L$ to
the empty link. Noting that $\sbn iL \cong \sbn i{\ol L}$ and that
$\sbn *{L\sqcup \ol L} \cong \sbn *L \ot \sbn *{\ol L}$ we see this
induces a map $\langle - , - \rangle \colon \sbn *L \ot \sbn *L \ra
\FF$. As is familiar in the discussion of 1+1-dimensional TQFTs the
diffeomorphism shown below implies that this inner product is
non-degenerate.

\begin{figure}[h]
\centerline{
\includegraphics{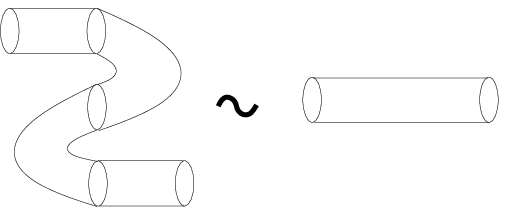}
}
\end{figure}

 This induces a
non-degenerate inner product on the complex $C^*$. 
Let $d_i^*$ denote the adjoint of $d_i$ with respect to this inner
product.

\begin{prop}
$ \obd_i = d^*_{i-1}$.
\end{prop}

\begin{proof}
In any 2d TQFT the linear map induced by a cylinder commutes with the linear 
map induced by any cobordism. Therefore the adjoint 
of the linear map induced by an oriented cobordism is given by the linear 
map induced by the same cobordism with the opposite orientation.   
\end{proof}

\begin{thm}\label{thm:conj}
\[
H_i(C^*, \obd)  \cong H_i(C^*, d)
\]
\end{thm}

\begin{proof}
By linear algebra (analogous to the discussion of harmonic
forms in Hodge theory)  we 
see that
\[
H_i(C^*, d) \cong \Ker (d_i) \cap \Ker (d^*_{i-1})
\]
and similarly
\[
H_i(C^*, \obd) \cong \Ker (\obd_i) \cap \Ker (\obd^*_{i+1}).
\]
However using the Proposition above we see that this becomes
\[
H_i(C^*, \obd) \cong \Ker (d^*_{i-1}) \cap \Ker ((d^*_{i})^*) =
\Ker (d^*_{i-1}) \cap \Ker (d_{i}) \cong H_{i}(C^*, d).
\]
\end{proof}

\vspace*{1cm}

\noindent {\bf Acknowledgements} The first author was supported by the 
Funda\c {c}\~{a}o para a Ci\^{e}ncia e a Tecnologia through the programme ``Programa 
Operacional Ci\^{e}ncia, Tecnologia, Inova\c {c}\~{a}o'' (POCTI), cofinanced by 
the European Community fund FEDER. 

The second author was supported by the
European Commission through a Marie Curie fellowship and thanks the
Institut de Recherche Math\'ematiques Avanc\'ee in Strasbourg for
their hospitality. 


\end{document}